\documentclass[12pt,reqno,a4paper]{amsart}
\usepackage{verbatim,rcs,cite}
\usepackage{amssymb,amsfonts,latexsym,mathrsfs}
\usepackage{amsmath}
\usepackage{amsthm}
\usepackage{mathrsfs}
\usepackage[all]{xy}
\usepackage{pstricks}

\theoremstyle{remark}

\theoremstyle{definition}

\DeclareMathOperator{\Sym}{Sym}

\RCS $Revision:  $ \RCSdate $Date: 2008/01/29 13:25:27 $

\title[ Magma Proof of Minimal Degrees]{Magma Proof of Strict Inequalities for Minimal Degrees of Finite Groups}

\date{}

\author{Scott H. Murray}
\author{Neil Saunders}

\dedicatory{\upshape
School of Mathematics and Statistics \\
University of Sydney, NSW 2006, Australia\\[.5em]
{\it E-mail address:} \ \texttt{S.Murray@maths.usyd.edu.au \\ neils@maths.usyd.edu.au }\\[1em] 
}

\begin{document}



\maketitle

\begin{abstract}
The minimal faithful permutation degree of a finite group $G$, denote by $\mu(G)$ is the least non-negative integer $n$ such that $G$ embeds
inside the symmetric group $\Sym(n)$. In this paper, we outline a Magma proof that $10$ is the smallest degree for which there are groups $G$ and $H$ such that $\mu(G \times H) < \mu(G)+ \mu(H)$.
\end{abstract}

\section{Introduction}
The study of this topic dates back to Johnson \cite{J71} and Wright \cite{W75}, who among other things investigated the
inequality
\begin{center}
\begin{equation} \label{eq:directsum}
\mu(G\times H) \leq \mu(G) + \mu(H),
\end{equation}
\end{center}
which clearly holds.\vspace{12pt}

Johnson first showed that equality holds when $G$ and $H$ have coprime orders or are abelian.  Wright went further to show that equality holds when $G$ and $H$ are $p$-groups and hence extended this to nilpotent groups. In that same paper, he constructs a class of groups $\mathscr{C}$ with the defining property that for every $G$ in $\mathscr{C}$, there exists a nilpotent subgroup $G_1$ in $G$ such that $\mu(G_1)=\mu(G)$. It is clear that equality in  \eqref{eq:directsum} holds for any two groups in $\mathscr{C}$ and that $\mathscr{C}$ is closed under taking direct products.\vspace{12pt}

Wright \cite{W75}, asked the question: does $\mu(G \times H)= \mu(G) + \mu(H)$ for all finite groups $G$ and $H$? The referee to \cite{W75} provided an example of strict inequality of degree $15$ and attached it as an addendum to that paper. The second author of this article recognised that the example quoted in that paper involved the complex reflection group $G(5,5,3)$ and its centraliser in $\Sym(15)$.\vspace{12pt}

This led to the investigation in \cite{S07} where the second author proved that a similar result occurs with the complex reflection groups $G(4,4,3)$ and $G(2,2,5)$, which are of degree $12$ and $10$ respectively. That is, these groups have non-trivial centralisers in their minimal embedding that intersect trivially their embedded image. \vspace{12pt}

In \cite{S08}, the second author extended this idea exhibiting that for $p$ and $q$ distinct odd primes, with $q \geq 5$ or $q=3$ and $p \not\equiv 1$ mod $3$, the groups $G(p,p,q)$ and their centralisers in $\Sym(pq)$ have the same property that $$\mu(G(p,p,q))=\mu(G(p,p,q)\times C_{\Sym(pq)}(G(p,p,q)))=pq,$$ and so give examples of strict inequality in \eqref{eq:directsum}. \vspace{12pt}

The authors do not know whether there are groups $G$ and $H$ such that $$\max \{\mu(G),\mu(H)\} < \mu(G \times H) < \mu(G) + \mu(H).$$

In the following section, we prove using the computation algebra system Magma \cite{CB06}, that $10$ is the smallest degree for the scenario that $\mu(G)=\mu(G \times C)$ where $G$ is minimally embedded group in $\Sym(\mu(G))$ and $C$ is its centraliser which intersects trivially with it. This is done by a brute-force search of the subgroups of $\Sym(m)$ for $m \leq 9$ and examining their centralisers.

\section{The Magma Code}
The following code was implemented in magma for $m \leq 9$
\begin{verbatim}

n:=m;
S:=Sym(m); 
num:=NumberOfTransitiveGroups(m);
subs:=Subgroups(Sym(m));
subs:=[s`subgroup: s in Subgroups(Sym(m))];
smaller:=[[s`subgroup: s in Subgroups(Sym(i))] : i in [1..m-1]];
minemb:=[ G : G in subs | forall{H : H in smaller[i], 
i in [1..m-1]| not IsIsomorphic(G,H)}];
Ind:=[Index(sub<S|Centraliser(S,G),G>,G) : G in minemb];
indices_min:=[i : i in [1..#minemb]| Ind[i] ne 1];
\end{verbatim}

\vspace{12pt}
Thus the code constructs the entire subgroup lattice of the symmetric group, isolates the subgroups which are minimally embedded inside the symmetric group and then computes their centralisers. For $G$ a minimally embedded group in $\Sym(m)$ and $C$ the centraliser of $G$ in this minimal embedding, the Ind sequence returns the index of $G$ in the group generated $G$ and $C$. Once this index is known, one can either determine that the centraliser is contained inside the group, or there is a possibility of a subgroup in $C$ which intersects trivially with $G$ by searching for an element in $x$ in $C$, such that the intersection of $\langle x \rangle$ with $G$ is trivial.
\section{Results}

Since in the cases $m=2,3,4$ are easily dealt with by hand we only give the Magma output for the higher cases. \vspace{12pt}

\vspace{12pt}
 
\underline{$m=5$}

\begin{verbatim}
> Ind;
[ 1, 1, 1, 1, 1, 1, 1 ]



\end{verbatim}

\underline{$m=6$}

\begin{verbatim}
> Ind;
[ 1, 1, 1, 1, 1, 1, 1, 1, 1, 1, 1, 1, 1, 1, 1, 1, 1, 1 ]

\end{verbatim}

\underline{$m=7$}

\begin{verbatim}

> Ind;
[ 1, 1, 1, 1, 1, 2, 1, 1, 1, 1, 1, 1, 1, 1, 1, 1, 1, 1, 1, 1,
 1, 1, 1, 1, 1, 1, 1, 1, 1 ]
\end{verbatim}

\underline{$m=8$}

\begin{verbatim}
> Ind;
[ 1, 4, 1, 1, 1, 1, 1, 1, 1, 1, 1, 1, 1, 1, 1, 2, 1, 2, 1,
 2, 1, 1, 1, 1, 1, 1, 1, 1, 1, 1, 1, 1, 1, 1, 1, 1, 1, 1, 
1, 1, 1, 1, 1, 1, 1, 1, 1, 1, 1, 1, 1, 1, 2, 1, 1, 1, 1, 1,
 1, 1, 1, 1, 1, 1, 1, 1, 1, 1, 1, 1, 1, 1, 1, 1, 1, 1, 1, 
1, 1, 1, 1, 1, 1, 1, 1, 1, 1, 1, 1, 1, 1, 1, 1, 1, 1, 1, 1,
 1, 1, 1, 1, 1, 1,1, 1, 1, 1 ]

\end{verbatim}

\underline{$m=9$}
\begin{verbatim}
Ind;
[ 1, 1, 1, 1, 1, 1, 2, 1, 1, 1, 1, 1, 1, 1, 2, 1, 1, 1, 1, 
1, 1, 1, 1, 1, 1, 1, 1, 1, 1, 1, 1, 1, 1, 1, 1, 1, 1, 1, 
1, 1, 1, 1, 1, 1, 1, 1, 1, 1, 1, 1, 1, 1, 1, 1, 1, 1, 1, 
1, 1, 1, 1, 1, 1, 1, 1, 1, 1, 1, 1, 1, 1, 1, 1, 1, 1, 1, 
1, 1, 1, 1, 1, 1, 1, 1, 1, 1, 1, 1, 1, 2, 1, 1, 1, 1, 1, 
1, 1, 1, 1, 1, 1, 1, 1, 1, 1, 1, 1, 1, 1, 1, 1, 1, 1, 1, 1
, 1, 1, 1, 1, 1, 1, 1, 1, 1, 1, 1, 1, 1, 1 ]
\end{verbatim}
\vspace{12pt}

An inspection of the numbers above shows that they are either $1$ or divisible by $2$. This means that any subgroup of the centraliser of $G$ in $\Sym(n)$ which intersects trivially with $G$ must have order divisible by $2$. Therefore, to search for such a subgroup, we implement the following function;\vspace{12pt}

\begin{verbatim}
Comp:= [ G : G in minemb |exists{g : g in Centraliser(S,G) |
 Order(g) eq 2 and Order(G meet sub<S|g>) eq 1} ];
\end{verbatim}
\vspace{12pt}

In each case, we find that 

\begin{verbatim}
> Comp;
[]
\end{verbatim}

\vspace{12pt}

Thus for every minimally embedded group of degree at most $9$, there does not exist a subgroup of its centraliser which intersects it trivially. Therefore we cannot obtain a strict inequality in \eqref{eq:directsum} by this method.


\bibliographystyle{plain}

\vspace{12pt}

\end{document}